\documentclass[12pt]{article}
\usepackage{amssymb}
\usepackage{amsmath}
\def\R{{\mathbb{R}}}
\title{Integral Functionals of Probability Measures \\
that Depend Only on the Mean
}
\author{Daniel W.~Stroock}

\begin{document}

\maketitle

\noindent
In this note I prove the following statement:
\medskip

\noindent
\textit{Assume that $f:\R\longrightarrow \R$ is a continuous function with the
property that there is a function $F$ such that
$$\int f\,d\mu =F\left(\int y\,\mu (dy)\right)$$
for all compactly supported probability measures $\mu $ on $\R$.
Then $f$ is an affine function. }
\medskip

\noindent
\textbf{\textit{Proof:\/}} Let $\rho $ be any non-negative, $C^\infty $, compactly
supported, even function with Lebesgue integral $1$,
and set $\rho_\epsilon (x)=\epsilon ^{-1}\rho(\epsilon ^{-1}x)$ for $\epsilon >0$.  Then
$$\int\rho _\epsilon (x-y)\,dy=x\text{ for all  } \epsilon >0\text{ and } x\in\R,$$
and so $\rho _\epsilon *f=F$ for all $\epsilon >0$.  Hence $f=F\in C^\infty $.  Next,
$$\frac{f(a+h)+f(a-h)}2=\int f\,d\left(\frac{\delta _{a+h}+\delta _{a-h}}2\right)=
f(a)\text{ for all } a\in\R\text{ and } h>0.$$
Hence
$$f''(a)=\lim_{h\searrow0}\frac{f(a+h)+f(a-h)-2f(a)}{h^2}=0\text{ for all } a\in\R.$$
\end{document}